\documentclass[a4paper,12pt]{article}
\usepackage{graphics} 
\usepackage{epsfig} 

\usepackage{amsmath,amssymb}
\usepackage{amsthm}

\usepackage{enumerate}

\allowdisplaybreaks[1]

\newtheorem{thm}{Theorem}

\newtheorem{rem}{\rm REMARK}[section]
\newtheorem{df}{\rm DEFINITION}
\newtheorem{ex}{\rm EXAMPLE}
\newtheorem{prop}{Proposition}[section]

\newcommand{\III}{{\vert \kern-.10em \vert \kern-.10em \vert}}

\makeatletter
 
 \@addtoreset{equation}{section}
\makeatother

\title{\Large EFFECTS OF RANDOMIZATION ON ASYMPTOTIC PERIODICITY
OF NONSINGULAR TRANSFORMATIONS}
\author{HIROSHI ISHITANI and KENSUKE ISHITANI}
\date{}

\begin{document}
\maketitle

\begin{abstract}
\noindent
It is known that the Perron--Frobenius operators of
piecewise expanding $\mathcal{C}^2$ transformations
possess an asymptotic periodicity of densities.
On the other hand, external noise or
measurement errors are unavoidable in practical systems;
therefore, all realistic mathematical models should be regarded
as random iterations of transformations.
This paper aims to discuss the effects of randomization
on the asymptotic periodicity of densities.
\end{abstract}

\section{Introduction}\label{Sec_Intro}
It is known that if $T: [0, 1) \rightarrow [0, 1)$ is a piecewise expanding $\mathcal{C}^2$ transformation,
then its corresponding Perron--Frobenius operator $\mathcal{L}_T$,
which we define in Section \ref{sec_preliminaries},
exhibits an asymptotic periodicity of densities.
That is, there exist probability density functions $g_{i,j,T}\in L^1([0,1))$ and functionals
$\lambda_{i,j,T}(\cdot)$ on $L^{1}([0,1))$ $(1\leq i \leq s(T)$, $1\leq j \leq r(i,T))$
satisfying the following conditions (see Definition \ref{def_of_asymPeriodicity_for_L} in Section \ref{sec_preliminaries}):
\begin{enumerate}[{(}i{)}]
\item
$g_{i,j,T}\cdot g_{k,l,T}=0$ for all $(i, j)\neq(k, l)$;
\item
For every $i$ $\in $ $\{1, 2, \ldots , s(T)\}$, $\{g_{i,j,T}\}_{j=1}^{r(i,T)}$ are periodic points of $\mathcal{L}_T$: \\
$\mathcal{L}_T(g_{i,j,T}) = g_{i,j+1,T}$ $(1\leq j \leq r(i,T)-1)$
and $\mathcal{L}_T(g_{i, r(i,T),T})$ $=$ $g_{i,1,T}$ hold;
\item
$\displaystyle \lim_{n\rightarrow \infty }\big \Vert
(\mathcal{L}_{T})^n
\big( f -\sum_{i=1}^{s(T)} \sum_{j=1}^{r(i,T)} \lambda_{i,j,T} (f) g_{i,j,T}\big) \big
\Vert_{L^{1}([0,1))} =0$ for $f \in L^{1}([0,1))$.
\end{enumerate}
Recall that the asymptotic periodicity of $\mathcal{L}_T$
describes the ergodic properties of the transformation $T$ (see, for example, \cite{InI}).
If we define $A_{i,j}=\{ x\in [0,1) ; g_{i,j,T}(x) >0 \}$ and $A_i = \bigcup_{j=1}^{r(i,T)} A_{i,j}$,
then the asymptotic periodicity of densities for $\mathcal{L}_T$ [(i)--(iii)] implies that
$T$ exhibits the following asymptotic periodicity \cite{InI},\cite{LiY},\cite{Wa}:
\begin{enumerate}[{(}a{)}]
\item
For every $i$ $\in $ $\{1, 2, \ldots , s(T)\}$, $A_i$ is $T$-invariant [{\it i.e.}, $T(A_i) =A_i$],
and the restriction $T_{i}\equiv T|_{A_i}$ is ergodic with respect to the Lebesgue measure $m$;
\item
For every $i$ $\in $ $\{1, 2, \ldots , s(T)\}$, there exists a $T$-invariant measure $\mu_i$ on $A_i$
that is equivalent to $m\vert _{A_i}$;
\item
For $B\equiv [0, 1) \setminus \bigcup_{i=1}^{s(T)} A_i$, we have that
$T^{-1} (B) \subset B$, and $\lim_{n\rightarrow \infty }m(T^{-n}(B)) = 0$;
\item
For every $i$ $\in $ $\{1, 2, \ldots , s(T)\}$, we have that
for the power $T_i^* \equiv {(T_i)}^{r(i,T)}$ of the transformation $T_i$,
$T_i^*(A_{i, j})=A_{i, j}$ $(1\leq j\leq r(i,T))$ holds, and $T_i^*$
is an exact endomorphism on $A_{i,j}$ $(1\leq j \leq r(i,T))$.
Further, for every $i$ $\in $ $\{1, 2, \ldots , s(T)\}$,
the transformation $T_i$ permutes $\{A_{i, j}\}_{j=1}^{r(i,T)}$ cyclically:
$T_i (A_{i, j}) = A_{i, j+1}$ $(1\leq j \leq r(i,T)-1)$ and
$T_i (A_{i, r(i,T)}) = A_{i, 1}$ hold.
\end{enumerate}
The above argument outlines the asymptotic periodicity of a single transformation $T$.
On the other hand, external noise, measurement errors, or inaccuracy
are unavoidable in practical systems. Therefore,
every realistic mathematical model should be regarded
as a number of random iterations of transformations $T_y\ (y \in Y)$:
$$
T_{\omega _{n}}\circ T_{\omega _{n-1}}\circ \cdots \circ T_{\omega _{1}}x,
$$
which is the first coordinate of the iterations of the skew product transformation
\begin{equation*}
S^{n}(x,\omega )
=(T_{\omega _{n}}\circ T_{\omega _{n-1}}\circ \cdots \circ
T_{\omega _{1}}x, \sigma ^{n}\omega )
\end{equation*}
defined in Subsection \ref{subsec_randomiterations} \cite{Kaku}.
In this paper, we consider only the case in which
transformations $T_{\omega_i}$ are independently chosen.
Then, under some assumptions, the skew product transformation $S(\equiv S^{1})$
is known to have asymptotic periodicity in the sense given above.
In this regard, it must be noted that the skew product transformation $S$
can be regarded as a random transformation.
This paper is concerned with the effects of this type of randomization on the asymptotic periodicity.

This paper is organized as follows.
In Section \ref{sec_preliminaries}, we review the necessary concepts
and results from the general theory of Perron--Frobenius operators,
as well as those relating to the random iterations of nonsingular transformations.
In Section \ref{sec_MainResults}, our main results are presented.
In Section \ref{sec_SufficientCondition},
we discuss a sufficient condition for the assumption of our main results.
In Section \ref{sec_ExamplesAndApplications}, we present some examples
with numerical experiments.

\section{Preliminaries}\label{sec_preliminaries}

In Subsection \ref{subsec_PF_operator}, we define 
the Perron--Frobenius operator and state its basic properties
that are necessary for our discussion.
In Subsection \ref{subsec_randomiterations},
we review the necessary concepts and results
from the theory of random iterations of transformations.

Although most of the results in this section are already well known
or can be easily seen, we give some of their proofs for completeness.

\subsection{Perron--Frobenius operators}\label{subsec_PF_operator}

Let $(X, \mathcal{F}, m)$ be a probability space and
$T: X\rightarrow X$ be an $m$-nonsingular transformation, {\it i.e.},
a measurable transformation, satisfying $m(T^{-1}(A)) = 0$
for $A\in \mathcal{F}$ with $m(A) = 0$.
Further, we denote the set of $p$-th integrable functions on $X$
with respect to the measure $m$ as
$L^p(m)\equiv L^p(X, \mathcal{F}, m)$ $(p\in [1,\infty ])$.
Then, we define the Perron--Frobenius operator corresponding to
$(X, \mathcal{F}, m, T)$ as follows.
\begin{df}
The Perron--Frobenius operator $\mathcal{L}_{T}$ on $L^1(m)$ is defined as 
\begin{equation}
\mathcal{L}_{T}f \equiv \frac{dm_f}{dm}, \quad
\mbox{where}\ m_f(A) = \int _{T^{-1}(A)} f(x) dm(x).
\nonumber
\end{equation}
\end{df}

The Perron--Frobenius operator $\mathcal{L}_{T}:L^1(m)\rightarrow L^1(m)$
is characterized by the following well-known proposition:
\begin{prop}\label{prop:2.1}
For $f$ $\in L^1(m)$, $\mathcal{L}_{T}f$ is the unique element in $L^1(m)$ satisfying
\begin{equation*}
\int _X (\mathcal{L}_{T}f)(x) h(x) dm(x) = \int _X f(x) h(Tx) dm(x)
\end{equation*}
for every $h \in L^{\infty}(m)$.
\end{prop}

As an operator on $L^1(m)$, $\mathcal{L}_{T}$ has the following properties, 
which are easily shown from Proposition \ref{prop:2.1}:
\begin{prop}\label{proposition2.2_PFoperator}
The operator $\mathcal{L}_{T}$ on $ L^1(m)$ is positive,
bounded, and linear, and it has the following properties:
\begin{enumerate}[{(}1{)}]
\item
$\mathcal{L}_{T}$ preserves integrals; {\it i.e.},
$\displaystyle \int_X (\mathcal{L}_{T}f)(x) dm(x) = \int_X f(x) dm(x)$
holds for $f \in L^1(m)$;
\item
For $f\in L^1(m)$, we have the inequality
$|(\mathcal{L}_{T}f)(x)| \leq (\mathcal{L}_{T} |f|)(x)$ $(m\mbox{-a.e.})$;
\item
$\mathcal{L}_{T}$ is a contraction; {\it i.e.},
$\Vert \mathcal{L}_{T}f \Vert _{L^1(m)} \leq \Vert f \Vert _{L^1(m)}$
holds for $f \in L^1(m)$;
\item
$(\mathcal{L}_{T})^n = \mathcal{L}_{T^n}$ holds,
where $\mathcal{L}_{T^n}$ represents the
Perron--Frobenius operator corresponding to $T^n$;
\item
For $g\in L^{\infty }(m)$ and $f\in L^{1}(m)$, we have
$g \cdot \mathcal{L}_{T}f=\mathcal{L}_{T}((g\circ T)f)$,
where \\$(g\circ T)(x)\equiv g(Tx)$;
\item
$\mathcal{L}_{T}f = f$ if and only if $f(x) dm(x)$ is $T$-invariant.
\end{enumerate}
\end{prop}

By applying Proposition \ref{proposition2.2_PFoperator},
we can obtain Propositions \ref{prop:2} and \ref{prop:3}.

\begin{prop}\label{prop:2}
Assume that $\mathcal{L}_{T}f = g$ holds for some nonnegative functions
$f$, $g\in L^{1}(m)$. Then
$$\displaystyle T^{-1}\{ g > 0 \} \supset \{ f > 0 \}\quad (m\mbox{-a.e.});$$
{\it i.e.},
$$\displaystyle m\left(\{ f>0 \} \setminus T^{-1} \{ g > 0 \} \right) =0$$
is satisfied.
\end{prop}
\proof
Using Proposition \ref{proposition2.2_PFoperator} (5) and
the assumption that $\mathcal{L}_{T}f =g$, we can show that
\begin{equation*}
(\mathcal{L}_{T}f )(x)=g(x)=1_{\{ g>0\}}(x) g(x)
= 1_{\{ g>0\}} (x) (\mathcal{L}_{T}f)(x)
= \mathcal{L}_{T}(1_{T^{-1}\{ g>0\}}f)(x).
\end{equation*}
We have $\int_X f(x) dm(x) = \int_X 1_{T^{-1}\{ g>0\}}(x) f(x) dm(x)$,
as $\mathcal{L}_{T}$ preserves integrals.
Therefore, the inequality $f(x) - 1_{T^{-1}\{ g>0\}}(x)f(x)\geq 0$ $(x\in X)$
shows that $f(x)= 1_{T^{-1}\{ g>0\}}(x)f(x)$\ $(m$-a.e.$\ x)$.
This completes the proof.
\qed

\begin{prop}\label{prop:3}
Let $f\in L^{1}(m)$ be nonnegative, and let $A\in \mathcal{F}$.
If $\mathcal{L}_{T}f = f$ and $T^{-1}(A) \supset A$, then $\mathcal{L}_{T}(f 1_A) = f 1_A$.
\end{prop}
\proof
Using the given assumptions
and Proposition \ref{proposition2.2_PFoperator} (5),
we have that
\begin{equation}\label{ineq_2.4}
1_A (x) f(x) = 1_A (x) (\mathcal{L}_T f)(x)
= \mathcal{L}_T(f 1_{T^{-1}(A)})(x)
\geq \mathcal{L}_T(f 1_{A})(x),\quad (m\mbox{-a.e.}\ x).
\end{equation}
By combining inequality (\ref{ineq_2.4}) and the fact that
$\int_X \left\{ 1_A (x) f(x) -\mathcal{L}_T (1_A f)(x) \right\} dm(x) = 0$,
we obtain $\mathcal{L}_{T}(f 1_A) = f 1_A$.
\qed

If  $\lim_{n\rightarrow \infty} (\mathcal{L}_{T})^{n} f = g$ holds,
then the limit set of $T^n\{ f\neq 0 \}$ is the support of $g$.
That is, we have the following proposition.

\begin{prop}\label{prop:2.5}
Assume that
\begin{equation}\label{prop2.5_eq:1}
\lim _{n\rightarrow \infty}
\Vert (\mathcal{L}_{T})^{n} f - g \Vert_{L^1(m)} =0
\end{equation}
holds for some nonnegative functions $f$, $g\in L^{1}(m)$.
Then
\begin{equation}\label{prop2.5_eq:1-1}
m\Big( \{ f > 0 \}
\setminus \bigcup _{n=0}^{\infty} T^{-n}\{ g > 0 \} \Big)=0.
\end{equation}
\end{prop}
\proof
First, we prove that the equation
\begin{equation}\label{prop2.5_proof_eq_1}
\lim_{n\rightarrow \infty}
m \left( \{ f>\varepsilon \} \setminus T^{-n}\{ g>0\} \right)= 0
\end{equation}
holds for any $\varepsilon > 0$.
In fact, we clearly have the estimate
\begin{align}
\Vert (\mathcal{L}_{T})^{n} f -
(\mathcal{L}_{T})^{n} \left(f 1_{T^{-n}\{ g>0\}} \right) \Vert_{L^1(m)}
& = \int_X ((\mathcal{L}_{T})^{n} \left( f - f 1_{T^{-n}\{ g>0\}} \right))(x) dm(x) \nonumber \\
& = \int_X f(x) \left( 1 - 1_{T^{-n}\{ g>0\}}(x) \right) dm(x) \nonumber \\
&  \geq \varepsilon \int_{\{f>\varepsilon \}}
\left( 1 - 1_{T^{-n}\{ g>0\}}(x) \right) dm(x) \nonumber \\
& = \varepsilon \ m \left( \{ f>\varepsilon \} \setminus T^{-n}\{ g>0\} \right).
\label{prop2.5_proof_ineq_1}
\end{align}
On the other hand, the inequality
\begin{align*}
& \Vert (\mathcal{L}_{T})^n f -
(\mathcal{L}_{T})^n \left( f 1_{T^{-n}\{ g>0\}} \right) \Vert_{L^1(m)} \\
&\quad \leq \Vert  (\mathcal{L}_{T})^{n} f -  g \Vert_{L^1(m)} +
\Vert  (\mathcal{L}_{T})^n \left( f 1_{T^{-n}\{ g>0\}} \right) -
g 1_{\{ g>0\}}  \Vert_{L^1(m)} \\
&\quad = \Vert  (\mathcal{L}_{T})^n f -  g \Vert_{L^1(m)} +
\Vert 1_{\{ g>0\}} \cdot \left(  (\mathcal{L}_{T})^n f - g \right) \Vert_{L^1(m)}  \\
&\quad \leq 2 \Vert (\mathcal{L}_{T})^n f - g \Vert_{L^1(m)},
\end{align*}
together with assumption $(\ref{prop2.5_eq:1})$, shows that
\begin{equation}\label{prop2.5_proof_eq_2}
\lim_{n\rightarrow \infty} \Vert
(\mathcal{L}_{T})^{n} f - (\mathcal{L}_{T})^{n} \left(f 1_{ T^{-n}\{ g>0\} } \right)
\Vert_{L^1(m)} =0.
\end{equation}
Therefore, the convergence in (\ref{prop2.5_proof_eq_1}) follows from
($\ref{prop2.5_proof_ineq_1}$) and ($\ref{prop2.5_proof_eq_2}$).

By applying assumption ($\ref{prop2.5_eq:1}$), we have that $\mathcal{L}_{T} g = g$.
As a result, Proposition \ref{prop:2} shows that
$T^{-n} \{ g>0 \} \supset T^{-(n-1)} \{ g>0 \}$ $(m\mbox{-a.e.})$ for $n\geq 1$.
Therefore, equation (\ref{prop2.5_proof_eq_1}) implies
\begin{equation}
m \left( \{ f>\varepsilon \} \setminus
\bigcup_{n=0}^{\infty}T^{-n}\{ g>0\} \right) = 0
\end{equation}
for any $\varepsilon >0$. This proves our result (\ref{prop2.5_eq:1-1}).
\qed

Let us define the asymptotic periodicity of $\mathcal{L}_T$ as follows.
\begin{df}\label{def_of_asymPeriodicity_for_L}
$\mathcal{L}_T$ is asymptotically periodic if there exist positive integers
$s(T)$ and $r(i, T)$ $(1\leq i \leq s(T))$, probability density functions $g_{i,j,T}\in L^1(m)$,
and bounded linear functionals $\lambda_{i,j,T}(\cdot )$ on $L^1(m)$
$(1\leq i \leq s(T), 1\leq j\leq r(i, T))$ such that
\begin{enumerate}[{(}i{)}]
\item
$g_{i,j,T}\cdot g_{k,l,T}=0$ for all $(i, j)\neq(k, l)$;
\item
$\mathcal{L}_T(g_{i, j, T}) = g_{i, j+1, T}$ $(1\leq j \leq r(i, T)-1)$ and
$\mathcal{L}_T(g_{i, r(i, T), T})=g_{i,1, T}$ hold for all $i$;
\item
$\displaystyle \lim_{n\rightarrow \infty }
\big\Vert  (\mathcal{L}_{T})^n
\big( f -\sum_{i=1}^{s(T)} \sum_{j=1}^{r(i, T)} \lambda_{i, j, T} (f) g_{i,j,T}\big)
\big\Vert_{L^{1}(m)} =0$ for any $f \in L^{1}(m)$.
\end{enumerate}
\end{df}

Using the above propositions and assuming 
the asymptotic periodicity of $\mathcal{L}_T$,
we can show the asymptotic periodicity of limit sets for $T$.

\begin{prop}\label{prop:2.6(r)}
Suppose that $\mathcal{L}_T$ is asymptotically periodic.
Then, if we denote \\
$\displaystyle g_{i,T}(x) \equiv \frac{1}{r(i,T)}\sum _{j=1}^{r(i,T)}g_{i,j,T}(x)$ $(x\in X)$,
$A_i \equiv \{ x \in X; g_{i,T}(x)>0 \}$, and \\$A_{i,j} \equiv \{ x \in X; g_{i,j,T}(x)>0 \}$,
$T$ has the following asymptotic periodicity:
\begin{enumerate}[{(}a{)}]
\item
For every $i$ $\in$ $\{1, 2, \ldots , s(T)\}$, $A_i$ is $T$-invariant [{\it i.e.}, $T(A_i) =A_i$],
and \\$d\mu_i(x)\equiv g_{i,T}(x) dm(x)$ is an ergodic $T$-invariant probability measure on $A _i$;
\item
For $B\equiv X\setminus \bigcup_{i=1}^{s(T)}A_i$, we have that
$T^{-1} (B) \subset B$ and $\lim_{n\rightarrow \infty }m(T^{-n}(B)) = 0$;
\item
For every power $T_i^*$ $\equiv$ ${(T_i)}^{r(i,T)}$ of $T_i\equiv T\vert _{A_i}$,
where $i$ $\in$ $\{1, 2, \ldots , s(T)\}$, we have that
$T_i^*(A_{i, j})=A_{i, j}$ $(1\leq j \leq r(i,T))$ holds, and
$T_i^*$ is an exact endomorphism on $A_{i, j}$ $[1\leq j \leq r(i,T)]$.
Further, for every $i$ $\in$ $\{1, 2, \ldots , s(T)\}$,
$T_i$ permutes $\{A_{i, j}\}_{j=1}^{r(i,T)}$ cyclically:
$T_i (A_{i, j}) = A_{i, j+1}$ $(1\leq j \leq r(i,T)-1)$
and $T_i (A_{i, r(i,T)}) = A_{i, 1}$ hold.
\end{enumerate}
\end{prop}

The following key proposition is established on the basis of the ergodicity
of each $A_i$, where $\{ A_1, A_2, \cdots , A_{s(T)}, B\}$
is the disjoint decomposition of $X$ given in Proposition \ref{prop:2.6(r)}.

\begin{prop}\label{prop:2-6}
Let $\{ A_1, A_2, \cdots , A_{s(T)}, B \}$ be a measurable partition
given in Proposition $\ref{prop:2.6(r)}$,
and let $A\in \mathcal{F}$ with $m(A)>0$ and $T^{-1}(A) \supset A$. Then,
\begin{align*}
&A_i \cap A = \emptyset \quad \mbox{or}\quad A_i \subset A
\quad \mbox{for}\quad i= 1, 2, \ldots , s(T),
\quad \mbox{and}\quad A \cap \bigcup _{i=1}^{s(T)} A_i \neq \emptyset .
\end{align*}
\end{prop}
\proof
Recall that $g_{i,T}(x)$ is a probability density function of
an ergodic, $T$-invariant measure $\mu_i$ on $A_i$.
From Proposition $\ref{prop:3}$, we obtain that $\mathcal{L}_T (g_{i,T} \cdot 1_A)=g_{i,T} \cdot 1_A$.
The ergodicity of $\mu_i$ allows us to state 
that either $g_{i,T} \cdot 1_A=0$ or $g_{i,T} \cdot 1_A = g_{i,T}$ holds for $1\leq i \leq s(T)$.
Therefore, either $A_i \cap A = \emptyset$ or $A_i \subset A$
holds for $1\leq i \leq s(T)$.
If we assume that $A \cap \bigcup _{i=1}^{s(T)} A_i = \emptyset$ holds,
then we have that $A \subset B$.
Under this assumption, we have that $T^{-n}(A) \supset A$
for $n$ $\in$ $\mathbb{N}$; hence,
$m(T^{-n}(B)) \geq m(T^{-n}(A)) \geq m(A) > 0$ for $n$ $\in$ $\mathbb{N}$.
This contradicts the condition that $\lim_{n\rightarrow \infty } m(T^{-n}(B))=0$.
\qed

\subsection{Random iteration}\label{subsec_randomiterations}

To formulate our main results, we need to introduce several further concepts.
In this subsection, we define the random iteration of $m$-nonsingular transformations.
\begin{enumerate}[{[}I{]}]
\item
Let $Y$ be a complete separable metric space, ${\cal B}(Y)$
be its topological Borel field, and $\eta $
be a probability measure on ($Y, {\cal B}(Y)$).
Further, define $\Omega \equiv \Pi _{i=1}^{\infty }Y$,
and let us write ${\cal B}(\Omega )$
for the topological Borel field of $\Omega $.
We insert the product measure
$P\equiv \Pi _{i=1}^{\infty }\eta $ on $(\Omega , {\cal B}(\Omega ))$.
\item
Let $(X, \mathcal{F}, m)$ be a probability space and $(T_{y})_{y\in Y}$
be a family of $m$-nonsingular transformations on $X$, such that
the mapping $(x,y) \rightarrow T_{y}x$ is measurable.
\end{enumerate}
To study the behavior of the random iterations,
we consider the skew product transformation
$S:X\times \Omega  \rightarrow  X\times \Omega $, defined as
\begin{equation}
S(x, \omega )\equiv (T_{\omega _{1}}x, \sigma \omega ),
\quad (x, \omega ) \in X \times \Omega ,
\end{equation}
where $\omega _{1}$ represents the first coordinate of
$\omega =(\omega _i)_{i=1}^{\infty }$, and
$\sigma : \Omega \rightarrow \Omega $
is the shift transformation to the left, which is defined as
$\sigma ((\omega _i)_{i=1}^{\infty })=(\omega _{i+1})_{i=1}^{\infty }$.
Note that, for $n \in \mathbb{N}$, we have
\begin{equation}
S^{n}(x,\omega )
=(T_{\omega _{n}}\circ T_{\omega _{n-1}}\circ \cdots \circ
T_{\omega _{1}}x, \sigma ^{n}\omega ).
\end{equation}
Therefore, we can consider the random iteration of $m$-nonsingular transformations as
$\pi S^{n}(x,\omega )$, writing $\pi : X\times \Omega \rightarrow  X$
for the projection on $X$.
Under these settings, Morita \cite{M1},\cite{M2},\cite{M3} investigated
the existence of invariant measures and their mixing properties.
His method is also useful for our purpose.

Because $(T_{y})_{y\in Y}$ are $m$-nonsingular transformations,
$S$ is a nonsingular transformation on
$(X \times \Omega , \mathcal{F}\times {\cal B}(\Omega ), m\times P)$.
Therefore, we can define the Perron--Frobenius operator
${\cal L}_S: L^{1}(m\times P) \rightarrow L^{1}(m\times P)$
corresponding to $S$ as
\begin{align*}
&\int \int _{X \times \Omega}
h(x, \omega ) ({\cal L}_S f) (x, \omega ) dm(x) dP(\omega )
= \int \int _{X \times \Omega}
f(x,\omega )h(S(x,\omega ))dm(x) dP(\omega ) \\
&\quad \mbox{for } h \in L^{\infty }(m\times P),
\mbox{ where }L^p(m\times P) \equiv
L^p(X \times \Omega , \mathcal{F}\times {\cal B}(\Omega ), m\times P)
\mbox{ for } p \in [1, \infty].
\end{align*}

Lemma 4.1 in \cite{M3} can be rewritten as follows:
\begin{prop}\label{prop:2.7}
\begin{enumerate}[{(}i{)}]
\item
If $({\cal L}_S f)(x,\omega ) = \lambda f(x,\omega )$ holds for $| \lambda | =1$,
then $f$ does not depend on $\omega $.
\item
For every $f \in L^{1}(m)$, we have
\begin{equation}
({\cal L}_S f)(x,\omega ) =
\int_Y ({\cal L}_{T_y} f)(x)\eta(dy), \quad (m\times P\mbox{-a.e.});
\end{equation}
hence, ${\cal L}_S f \in L^{1}(m)$.
\end{enumerate}
\end{prop}

Proposition \ref{prop:2.7} allows us to consider ${\cal L}_S$ as an operator on $L^{1}(m)$.
Then, we have the following key proposition:
\begin{prop}\label{prop:2.8}
Assume that $(\mathcal{L}_S f)(x,\omega) = g(x)$ $(m\times P\textrm{-$a.e.$})$
holds for some nonnegative functions $f$, $g\in L^{1}(m)$.
Then, there exists a set $Y_0 \in {\cal B}(Y)$, with $\eta (Y_0)=1$,
such that for every $\omega_1 \in Y_0$,
\begin{align*}\label{eq:}
&
(T_{\omega_1})^{-1}\{ g > 0 \}
\supset \{ f > 0 \} \ (m\mbox{-a.e.});
\quad
\mbox{{\it i.e.},}\quad m\left( \{ f > 0 \} \setminus
(T_{\omega_1})^{-1}\{ g > 0 \} \right) =0
\end{align*}
is satisfied.
\end{prop}
\proof
By applying Proposition \ref{prop:2}, we obtain the equation
\begin{equation*}
(m\times P)\left( \{ (x,\omega) \in X\times \Omega ;  \,f(x) > 0 \}
\setminus S^{-1}\{ (x,\omega) \in X\times \Omega ;  \,g(x) > 0 \} \right) =0.
\end{equation*}
Fubini's theorem implies that there exists a set $\Omega_0 \in {\cal B}(\Omega )$,
with $P(\Omega_0)=1$, such that
\begin{equation*}
m\left( \{ x\in X ; f(x) > 0 \} \setminus
(T_{\omega_1})^{-1}\{ x\in X ; g(x) > 0 \} \right) =0
\end{equation*}
holds for $(\omega_{i})_{i=1}^{\infty } \in \Omega_0$.
Define $Y_0 \equiv \widetilde{\pi}_1(\Omega_0)$, where $\widetilde{\pi}_1:\Omega \rightarrow Y$
is the projection given by $\widetilde{\pi}_1((\omega_{i})_{i=1}^{\infty })=\omega_1$
for $(\omega_{i})_{i=1}^{\infty } \in \Omega$.
Then we obtain our proposition.
\qed

\section{Main results}\label{sec_MainResults}

In this section, we state our main results using the notation
defined in Subsection \ref{subsec_randomiterations}.
To state our main results, we assume the asymptotic periodicity of $\mathcal{L}_S$.
Then there exist probability density functions $g_{i,j,S}$ $\in $ $L^1(m)$
and functionals $\lambda_{i,j,S}(\cdot)$ on $L^1(m)$
$(1\leq i \leq s(S)$, $1\leq j \leq r(i,S))$
satisfying the conditions of Definition \ref{def_of_asymPeriodicity_for_L}.
We denote $s(S)$, $r(i,S)$, $g_{i,j,S}$, and $\lambda_{i,j,S}(\cdot)$
by $\widehat{s}$, $\widehat{r}(i)$, $\widehat{g}_{i,j}$, and $\widehat{\lambda}_{i,j}(\cdot)$, respectively.
Note that $\widehat{s}$ is the number of ergodic components of $\mathcal{L}_{S}$,
and $\widehat{r}(i)$ is the number of cycles of each ergodic component $(1\leq i \leq \widehat{s})$.
It is also meaningful to consider a sufficient condition
for the above assumption,
which will be discussed in Section \ref{sec_SufficientCondition}.

Let $Y_1$ denote the set of parameters $y\in Y$ such that $\mathcal{L}_{T_y}$
is asymptotically periodic; that is,
for $y\in Y_1$, there exist probability density functions $g_{i,j,T_y}$ $\in L^1(m)$
and functionals $\lambda_{i,j,T_y}(\cdot)$ on $L^1(m)$
$(1\leq i \leq s(T_y)$, $1\leq j \leq r(i,T_y))$
satisfying the conditions of Definition \ref{def_of_asymPeriodicity_for_L}.
Note also that
$s(T_y)$ is the number of ergodic components of $\mathcal{L}_{T_y}$,
and $r(i,T_y)$ is the period of cycles of each ergodic component $(1\leq i \leq s(T_y))$.
If $\mathcal{L}_S$ is asymptotically periodic,
$\eta(Y_1)$ is positive in many cases; this is also confirmed by
the examples given in Section \ref{sec_ExamplesAndApplications}.
Under the above notation, Proposition \ref{prop:2.8} can be rewritten as follows:
\begin{prop}\label{prop:3-1}
Suppose that $\mathcal{L}_S$ is asymptotically periodic.
Then there exists a set $Y_0 \in {\cal B}(Y)$, with $\eta(Y_0)=1$, such that
\begin{align*}
&\{\widehat{g}_{i,j} > 0\} \subset T_y^{-1}\{\widehat{g}_{i, j+1} > 0\}
\ (1\leq j \leq \widehat{r}(i)-1)
\mbox{ and }
\{\widehat{g}_{i,\widehat{r}(i)} > 0\} \subset T_y^{-1}\{\widehat{g}_{i,1} > 0\}
\end{align*}
hold for $y\in Y_0$.
\end{prop}

\begin{rem}
Define
$\displaystyle \widehat{g}_i \equiv
\frac{1}{\widehat{r}(i)}\sum_{j=1}^{\widehat{r}(i)} \widehat{g}_{i,j}$
$(1\leq i \leq \widehat{s})$ and
$\displaystyle g_{i,T_y} \equiv \frac{1}{r(i,T_y)} \sum_{j=1}^{r(i,T_y)} g_{i,j,T_y}$
$(1\leq i \leq s(T_y))$.
Then, $\widehat{g}_i $ and $g_{i,T_y}$ are the densities of the ergodic invariant
probabilities of $S$ and $T_y$, respectively.
\end{rem}

We are now in a position to state the first main result of this paper.
\begin{thm}\label{thm:1}
Suppose that $\mathcal{L}_S$ is asymptotically periodic.
Let $Y_0$ be the set as in Proposition \ref{prop:3-1}.
Then, for every $i \in \{1, 2, \ldots , \widehat{s}\}$ and $y\in Y_0 \cap Y_1$,
we have the following statements.
\begin{enumerate}[{(}1{)}]
\item
We have that
$\{\widehat{g}_i>0\} \cap \bigcup _{k=1}^{s(T_y)} \{ g_{k,T_y}>0 \} \neq \emptyset$ holds.
Moreover, either $\{\widehat{g}_i>0\} \supset \{g_{k,T_y}>0\}$ or
$\{\widehat{g}_i>0\} \cap \{g_{k,T_y}>0\}=\emptyset$
holds for every $k$ $\in$ $\{1, 2, \ldots , s(T_y)\}$.
This means that $\widehat{s}$ is not greater than $s(T_y)$.
\item
For $k_0 \in \{1, 2, \ldots , s(T_y)\}$ satisfying
$\{ \widehat{g}_i>0 \} \supset \{ g_{k_0,T_y}>0 \}$,
$\widehat{r}(i)$ is a divisor of $r(k_0,T_y)$.
\end{enumerate}
\end{thm}
\proof
(1) Note that $\{ \widehat{g}_i > 0 \} = \bigcup _{j=1}^{\widehat{r}(i)} \{ \widehat{g}_{i,j}>0 \}$.
Then, Proposition \ref{prop:3-1} shows that
$T_y^{-1} \{ \widehat{g}_i > 0 \} \supset \{ \widehat{g}_i > 0 \}$.
Hence, we easily obtain the first statement, which follows from Proposition \ref{prop:2-6}.

\noindent
(2) For simplicity, let us write $\widehat{g}_{k,r}\equiv \widehat{g}_{k,j}$
for $r=l \cdot \widehat{r}(k)+j$
($1\leq k \leq \widehat{s}$, $1\leq j \leq \widehat{r}(k)$, $l \in \mathbb{N}$)
and $g_{k,r,T_y}\equiv g_{k,j,T_y}$
for $r=l \cdot r(k,T_y)+j$ ($1\leq k \leq s(T_y)$, $1\leq j \leq r(k,T_y)$, $l \in \mathbb{N}$).
Then, from Proposition \ref{prop:3-1},
we have that $\{ \widehat{g}_{k, j} > 0 \}$ $\subset$
$T_{y}^{-1}\{ \widehat{g}_{k, j+1} > 0 \}$ ($1\leq k \leq \widehat{s}$, $j\in \mathbb{N}$).

Because $g_{k_0,j+1,T_y}=\mathcal{L}_{T_y}(g_{k_0,j,T_y})$
holds for $j \in \mathbb{N}$, we have that
$$
1_{\{ \widehat{g}_{i, j+1}>0\}} g_{k_0,j+1,T_y}
=1_{\{ \widehat{g}_{i, j+1}>0\} }\mathcal{L}_{T_y}(g_{k_0,j,T_y})
=\mathcal{L}_{T_y}(1_{T_y^{-1}\{ \widehat{g}_{i, j+1}>0 \}}g_{k_0,j,T_y}).
$$
From the assumption that $\{ \widehat{g}_i>0 \} \supset \{ g_{k_0,T_y}>0 \}$,
we have that $g_{k_0,j,T_y}=1_{\{ \widehat{g}_{i} > 0\}}g_{k_0,j,T_y}$.
Note that the sets $\{ \widehat{g}_{i, j} > 0\}$ ($j=1, 2, \ldots, \widehat{r}(i)$)
are mutually disjoint. Then, we can obtain that
$T_{y}^{-1}\{ \widehat{g}_{i, j+1} > 0\}$ $\cap $ $\{ \widehat{g}_{i} > 0\}$
$=\{ \widehat{g}_{i, j} > 0\}$. Hence, we have that
$$
\mathcal{L}_{T_y}(1_{T_y^{-1}\{ \widehat{g}_{i, j+1}>0\}}g_{k_0,j,T_y})
=\mathcal{L}_{T_y}
(1_{T_y^{-1}\{ \widehat{g}_{i, j+1}>0\}}1_{\{ \widehat{g}_i > 0\}}g_{k_0,j,T_y})
=\mathcal{L}_{T_y}(1_{\{ \widehat{g}_{i, j}>0\}}g_{k_0,j,T_y}).
$$
This implies that
\begin{equation} \label{eq:3-1}
1_{\{ \widehat{g}_{i, j+1}>0\}}g_{k_0,j+1,T_y}
=\mathcal{L}_{T_y}(1_{\{ \widehat{g}_{i, j}>0\}}g_{k_0,j,T_y})
\quad \mbox{for}\ j \in \mathbb{N}.
\end{equation}
Recall the assumption that $\{ \widehat{g}_i>0 \} \supset \{ g_{k_0,T_y}>0 \}$.
Then, by renumbering, we can assume that
$1_{\{ \widehat{g}_{i, 1} > 0\}} g_{k_0,1,T_y}$ is a nontrivial function.
Thus, we obtain a sequence of nontrivial nonnegative functions
$h_j \equiv 1_{\{ \widehat{g}_{i, j} > 0 \}} g_{k_0,j,T_y}$ ($j \in \mathbb{N}$)
satisfying the equations $h_{j+1}=\mathcal{L}_{T_y}h_j$ ($j \in \mathbb{N}$)
and $h_j \cdot h_l = 0$ for $1\leq j <l \leq N$,
where $N$ denotes the least common multiple of $\widehat{r}(i)$ and $r(k_0,T_y)$.
It clearly follows that
\begin{align*}
\sum _{j=1}^{N} h_j &= \sum _{j=1}^{N} 1_{\{ \widehat{g}_{i, j} > 0 \}} g_{k_0,j,T_y}
\leq \sum_{r=1}^{\widehat{r}(i)} \sum_{j=1}^{r(k_0,T_y)}
1_{\{ \widehat{g}_{i, r} > 0 \}} g_{k_0,j,T_y}
=\sum_{j=1}^{r(k_0,T_y)} g_{k_0,j,T_y}.
\end{align*}
Therefore, we obtain the estimate
\begin{align*}
N\int_X h_1(x) dm(x) &=\sum_{j=0}^{N-1}\int_X ((\mathcal{L}_{T_y})^j h_1)(x) dm(x)\\
&= \int_X \sum _{j=1}^{N} h_j(x) dm(x) \\
&\leq \int_X \sum _{j=1}^{r(k_0,T_y)} g_{k_0,j,T_y}(x) dm(x) \\
&= r(k_0,T_y)\int_X g_{k_0,1,T_y}(x) dm(x).
\end{align*}
If $N > r(k_0,T_y)$ holds, then we have that
\begin{align*}
& \int_X h_1(x) dm(x)
\leq \frac{r(k_0,T_y)}{N}\int_X g_{k_0,1,T_y}(x) dm(x)
< \int_X g_{k_0,1,T_y}(x) dm(x).
\end{align*}
This shows that $\{h_1>0\} \subsetneqq \{g_{k_0,1,T_y}>0\}$; {\it i.e.},
$m(\{g_{k_0,1,T_y}>0\}\setminus \{h_1>0\})>0$.
Because $(\mathcal{L}_{T_y})^{n N}(h_1)=h_1$
and $(\mathcal{L}_{T_y})^{n N}(g_{i,j,T_y})=g_{i,j,T_y}$
hold for $n\in \mathbb{N}$, we have that
\begin{align*}
&\Big\Vert (\mathcal{L}_{T_y})^{n N}
\big( h_1 -
\sum_{k=1}^{s(T_y)} \sum_{j=1}^{r(k,T_y)} \lambda_{k,j,T_y}(h_1) g_{k,j,T_y}\big) \Big\Vert_{L^{1}(m)}
=\Big\Vert
h_1 -\sum_{k=1}^{s(T_y)} \sum_{j=1}^{r(k,T_y)} \lambda_{k,j,T_y}(h_1) g_{k,j,T_y} \Big\Vert_{L^{1}(m)} \\
&\quad = \Vert  h_1 - \lambda_{k_0,1,T_y}(h_1) g_{k_0,1,T_y} \Vert_{L^{1}(m)} +
\Big\Vert \sum_{(k, j)\neq (k_0,1)} \lambda_{k,j,T_y}(h_1) g_{k,j,T_y}
\Big\Vert_{L^{1}(m)} \\
&\quad \geq \Vert  h_1 - \lambda_{k_0,1,T_y}(h_1) g_{k_0,1,T_y}
\Vert_{L^{1}(m)} \\
&\quad = \int_X (1_{\{h_1>0 \}}(x) + 1_{\{h_1=0 \}}(x))
\big\vert h_1(x) - \lambda_{k_0,1,T_y}(h_1) g_{k_0,1,T_y}(x)\big\vert dm(x)\\
&\quad = \vert 1 - \lambda_{k_0,1,T_y}(h_1)
\vert \int_X h_1(x)\ dm(x) + \vert \lambda_{k_0,1,T_y}(h_1) \vert
\int_X 1_{\{g_{k_0,1,T_y}>0,\ h_1=0\}}(x) g_{k_0,1,T_y}(x) dm(x).
\end{align*}
If $\lambda_{k_0,1,T_y}(h_1)=1$ holds,
the right-hand side of the above inequality is
\begin{align*}
&\int_X 1_{\{g_{k_0,1,T_y}>0,\ h_1=0\}}(x) g_{k_0,1,T_y}(x) dm(x),
\end{align*}
which is strictly positive.
If $\lambda_{k_0,1,T_y}(h_1)\neq 1$ holds,
we have
\begin{align*}
&\vert 1 - \lambda_{k_0,1,T_y}(h_1) \vert \int_X h_1(x)\ dm(x)>0.
\end{align*}
Therefore, there exists a positive constant $a$ such that
\begin{align*}
&\Big\Vert
(\mathcal{L}_{T_y})^{nN}
\big(
h_1 -\sum_{k=1}^{s(T_y)} \sum_{j=1}^{r(k,T_y)} \lambda_{k,j,T_y}(h_1) g_{k,j,T_y}
\big) \Big\Vert_{L^{1}(m)}\geq a > 0
\end{align*}
holds for every $n \in \mathbb{N}$.
This contradicts the fact that
$$
\lim_{n\rightarrow \infty} \Big\Vert
(\mathcal{L}_{T_y})^{nN}
\big(
h_1 -\sum_{k=1}^{s(T_y)} \sum_{j=1}^{r(k,T_y)} \lambda_{k,j,T_y}(h_1) g_{k,j,T_y}
\big) \Big\Vert_{L^{1}(m)} = 0
$$
for $y \in Y_1$. Therefore, we have $N=r(k_0,T_y)$.
This implies that $\widehat{r}(i)$ is a divisor of $r(k_0,T_y)$.
\qed

When the identity map $I_d$ on $X$ is chosen with positive probability,
Proposition \ref{prop:3-1} can be applied to show that
the transformation $S$ is exact on $\{\widehat{g}_i>0\}$ $(1\leq i \leq \widehat{s})$:
\begin{thm}\label{thm:2}
Suppose that $\mathcal{L}_S$ is asymptotically periodic
and $\eta (\{y\in Y; \ T_y=I_d\})$ $>0$ is satisfied.
Then, for every $i\in \{1, 2, \ldots , \widehat{s}\}$,
$\widehat{r}(i)= 1$ holds.
\end{thm}
\proof
Assume that $\widehat{r}(i)\geq 2$ holds for some $1\leq i \leq \widehat{s}$.
Then, we have that $\mathcal{L}_{S}(\widehat{g}_{i, 1})=\widehat{g}_{i, 2}$ and
$\widehat{g}_{i, 1} \cdot \widehat{g}_{i, 2} = 0$.
It follows from Proposition \ref{prop:3-1} that
$\{ \widehat{g}_{i, 1}>0 \} \subset I_d^{-1}\{\widehat{g}_{i, 2}>0\}$.
This contradicts the fact that $\widehat{g}_{i, 1} \cdot \widehat{g}_{i, 2} = 0$.
\qed

\section{A sufficient condition for a skew product transformation
to have asymptotic periodicity of densities}\label{sec_SufficientCondition}

In this section, we discuss a sufficient condition for a skew product
of one-dimensional transformations to have asymptotic periodicity of densities.
Let $I\equiv [0,1]$ be the unit interval, $\mathcal{F}\equiv \mathcal{B}([0,1])$
be the Borel field, and $m$ be the Lebesgue measure on $(I, \mathcal{F})$.
We consider a family $(T_y)_{y\in Y}$ of $m$-nonsingular transformations on $[0,1]$,
where $Y$ is a complete separable metric space
equipped with a probability measure $\eta $ on $(Y, {\cal{B}}(Y))$.

For $f: [0,1]\rightarrow  \mathbb{C}$,
we denote the total variation of $f$ on $[0,1]$ by $var(f)$.
It is known that $V$ $\equiv$ $\{f \in L^{1}([0,1]); v(f)<\infty \}$
is a non-closed subspace of $L^{1}([0,1])$,
where $v(f)=\inf \{var(\tilde{f}); \tilde{f} \mbox{ is a version of } f\}$
for $f \in L^{1}([0,1])$.
On the other hand, letting
\begin{equation*}
\| f\| _{V}\equiv \| f\| _{L^1([0,1])} + v(f) \quad \mbox{for}\ f \in V,
\end{equation*}
we can easily prove that ($V$, $\| \cdot \| _{V})$
is a Banach space, and that the inequality \\
$\| f g\| _{V}\leq  2\| f\| _{V}\| g\| _{V}$ holds for $f, g \in V$ (cf. \cite{R}).

\begin{df}\label{def_D_infty}
Let ${\cal D}_{\infty }$ be the set of all transformations
$T:[0,1]\rightarrow [0,1]$ satisfying the following conditions:
\begin{enumerate}[{(}1{)}]
\item
There is a countable partition $\{ I_j \}_j $ of $I$ by disjoint intervals
such that the restriction $T\vert _{I_j}$ of $T$ to $I_j$
can be extended to a monotonic $\mathcal{C}^2$-function
on the closure $\bar{I_j}$ for each $j$,
and the collection $\{J_{j}\equiv T(I_j)\}_j$
consists of a finite number of different subintervals;
\item $\gamma (T)\equiv \inf_{x\in [0,1]} | T'(x) | > 0$ holds.
\end{enumerate}
\end{df}

We state the following inequality for a single transformation
$T$ $\in $ ${\cal D}_{\infty }$,
which was established by Rousseau-Egele $\cite{R}$.
\begin{prop}
Assume that $T\in {\cal D}_{\infty }$
and that the corresponding partition $\{I_j\}_j$
and $\gamma (T)$ from Definition \ref{def_D_infty}
are given.
Then, we have the following inequality:
\begin{align*}
&v({\cal L}_T f) \leq \alpha (T)v(f) + \beta (T)\| f\| _{L^1([0,1])},
\quad (f\in V), \\
&\quad \mbox{where}
\ \ \alpha (T)\equiv \frac{2}{\gamma (T)}
\ \ \mbox{and}
\ \ \beta (T)\equiv
\sup _j \Big\{ \frac{1}{m(I_{j})} \Big\} +
\sup_{j}\left\{
\frac{\sup_{x\in I_{j}}| (T_{j} ^{-1})''(x)|}{\inf_{x\in I_{j}}| (T_{j} )'(x)|}
\right\}.
\end{align*}
\end{prop}

We now consider the skew product transformation $S$
of $(T_y)_{y\in Y}$ $\subset $ ${\cal D}_{\infty }$.
The following proposition enables us to give a sufficient condition for ${\cal L}_S$ to
have asymptotic periodicity.

\begin{prop}
Assume that $(T_{y})_{y\in Y}\subset {\cal D}_{\infty }$ is given
and that the inequalities
\begin{align}
&\int_Y \int_Y \cdots \int_Y \alpha (T_{y_{n_0}}\circ T_{y_{n_0-1}}
\circ \cdots \circ T_{y_{1}})\eta(dy_{n_0})\eta(dy_{n_0-1})\cdots \eta(dy_{1}) < 1
\label{eq:alpha}\\
&\int_Y \int_Y \cdots \int_Y \beta (T_{y_{n_0}}
\circ T_{y_{n_0-1}}\circ \cdots \circ T_{y_{1}})
\eta(dy_{n_0})\eta(dy_{n_0-1})\cdots \eta(dy_{1}) < \infty
\label{eq:beta}
\end{align}
hold for some $n_0\in \mathbb{N}$.
Then, there exist real numbers $\alpha \in (0,1)$ and $\beta \in (0, \infty) $
satisfying
\begin{equation}\label{eq:qc}
v(({\cal L}_S)^{n_{0}}f) \leq \alpha  v(f) + \beta  \| f\| _{L^1([0,1])},
\quad (f \in V).
\end{equation}
\end{prop}

Using inequality \eqref{eq:qc}, the theorem of Ionescu-Tulcea and Marinescu \cite{IM}
showed the quasi-compactness and asymptotic periodicity of ${\cal L}_S$.
Note that the families of transformations in Section \ref{sec_ExamplesAndApplications}
satisfy inequalities \eqref{eq:alpha} and \eqref{eq:beta}; thus, 
$\mathcal{L}_S$ is asymptotically periodic.

\section{Numerical examples}\label{sec_ExamplesAndApplications}

This section uses examples to demonstrate our main results
from Section \ref{sec_MainResults}.
We consider the unit interval $X\equiv [0,1]$, Borel field
$\mathcal{F}\equiv \mathcal{B}([0,1])$,
and Lebesgue measure $m$ on $(X, \mathcal{F})$.
Further, in this section, we employ the initial density function
$f_0(x)=2x$ for $x\in [0,1]$,
the complete separable metric space $Y=\{y_1, y_2\}$
$\subset \mathbb{R}$\ ($y_1 \neq y_2$),
and the probability measure $\eta$ on $Y$ satisfying $\eta (\{y_1\})=\eta (\{y_2\})=1/2$.
Thus, the Perron--Frobenius operator ${\cal L}_S f$ is obtained as follows:
\begin{align*}
&(\mathcal{L}_S f)(x)=\frac{1}{2}\left\{
(\mathcal{L}_{T_{y_1}} f)(x)
+(\mathcal{L}_{T_{y_2}} f)(x)\right\}, \ x\in X.
\end{align*}
Because  $\eta (\{y_i\})>0$ $(i=1,2)$, it follows that $Y_0=\{y_1, y_2\}$.

\begin{ex}
For $m_0 \in \mathbb{N}$, we define the subintervals $J_k$
$(1$ $\leq $ $k$ $\leq $ $m_0)$ as follows:
\begin{align*}
&J_k \equiv \left[\frac{k-1}{m_0}, \frac{k}{m_0} \right), \ (1\leq k\leq m_0-1),
\ \mbox{and}\ J_{m_0}\equiv \left[1-\frac{1}{m_0}, 1\right].
\end{align*}
We consider the transformation $R_3$ on $X$ given by
$R_3 x \equiv 3x\ (\mbox{mod } 1)$.
Then, we define the transformation $R^{\tau}: X \rightarrow X$ as
\begin{align*}
&R^{\tau} x \equiv \frac{1}{m_0} R_3(m_0 x - k + 1)+\frac{\tau _k -1}{m_0},
\ \mbox{for}\ x \in J_k \ (k\in \{1, 2, \ldots , m_0\}),
\end{align*}
where $\tau = (\tau _1, \ldots , \tau _{m_0})$ is a permutation of the set
$\{1, 2, \ldots , m_0\}$.
Further, the Perron--Frobenius operator ${\cal L}_{R^{\tau}} f$ is obtained as follows:
\begin{align*}
({\cal L}_{R^{\tau}} f)(x) =& \frac{1}{3}\bigg\{
f\left(\frac{x}{3}+\frac{k-1}{m_0}-\frac{\tau_k -1}{3m_0}\right)
+ f\left(\frac{x}{3}+\frac{k-1}{m_0}-\frac{\tau_k -2}{3m_0} \right) \\
&\ \ + f\left(\frac{x}{3}+\frac{k-1}{m_0}-\frac{\tau_k -3}{3m_0} \right) \bigg\} 
\mbox{ for } x \in J_{\tau_k} = R^{\tau}(J_k)\ (1\leq k\leq m_0).
\end{align*}
The graph of $\{R^{\tau}x;\ x\in [0, 1]\}$ for $\tau = (2,1,6,4,3,5)$ is shown in
Fig \ref{fig:Graph_R_tau216435}.
\begin{figure}[htp]
\begin{center}
\rotatebox{-90}
{\includegraphics[width=43mm]{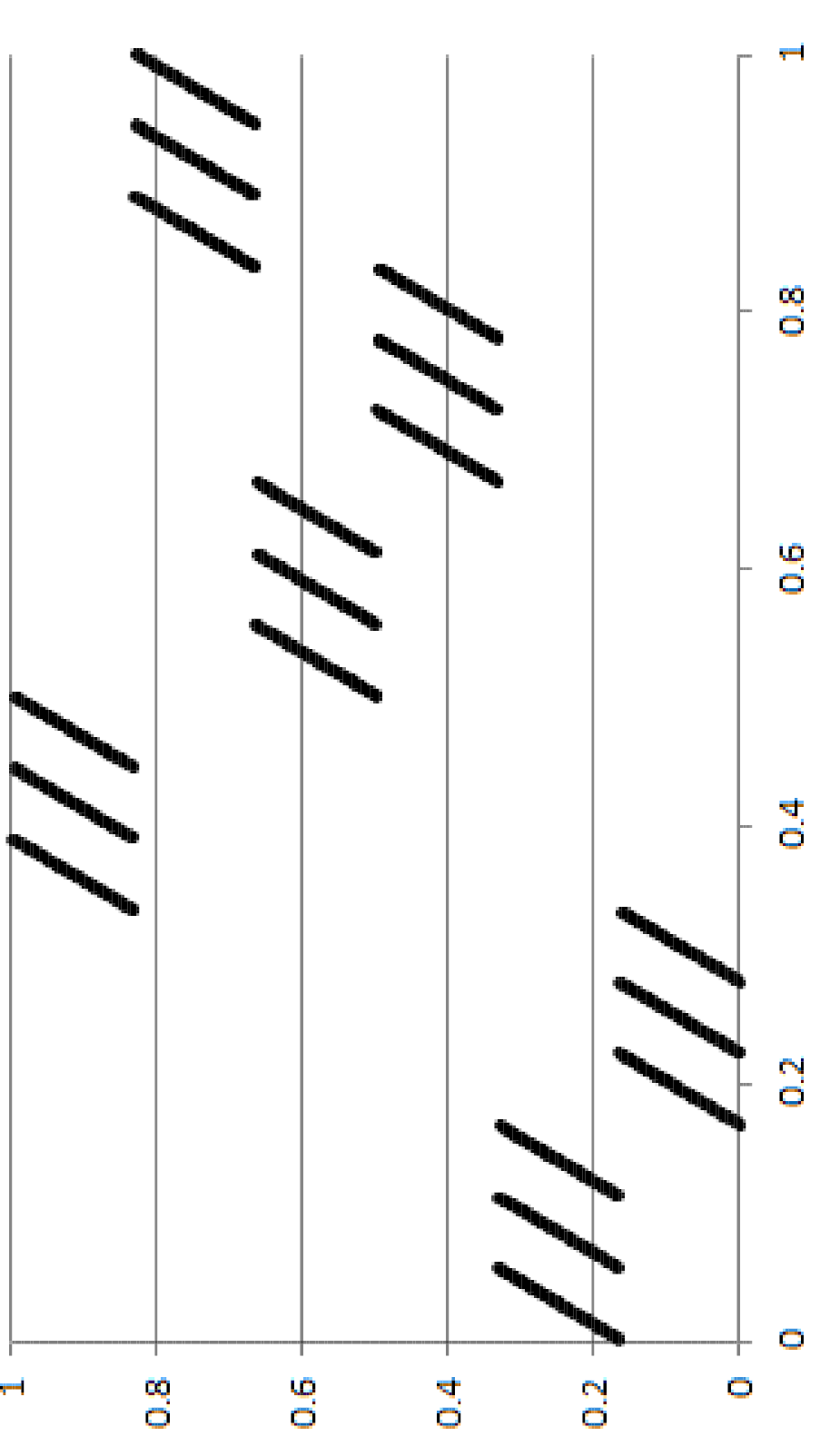}}
\end{center}
\caption
{$\{R^{\tau}x;\ x\in [0, 1]\}$ for $\tau = (2,1,6,4,3,5)$. }
\label{fig:Graph_R_tau216435}
\end{figure}

\subsubsection*{For $T_{y_1}=R^{(3,4,2,1)}$ and $T_{y_2}=R^{(4,5,6,2,3,1)}$}
\begin{figure}[htp]
\begin{center}
\rotatebox{-90}{\includegraphics[width=45mm]{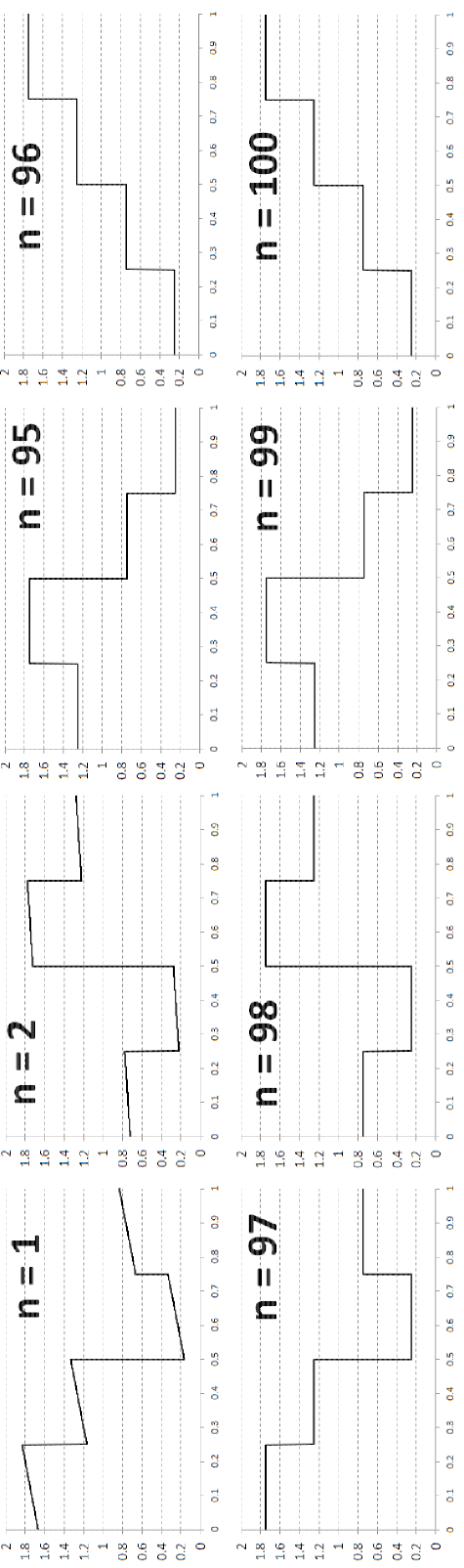}}
\end{center}
\caption
{Results of $((\mathcal{L}_{R^{\tau}})^n f_0)(\cdot )$
for $\tau $= $(3,4,2,1)$,
$n=1$, $2$, $95$, $96$, $97$, $98$, $99$, and $100$. }
\label{fig:R_tau3_3421}
\end{figure}
\begin{figure}[htp]
\begin{center}
\rotatebox{-90}{\includegraphics[width=37.9mm]{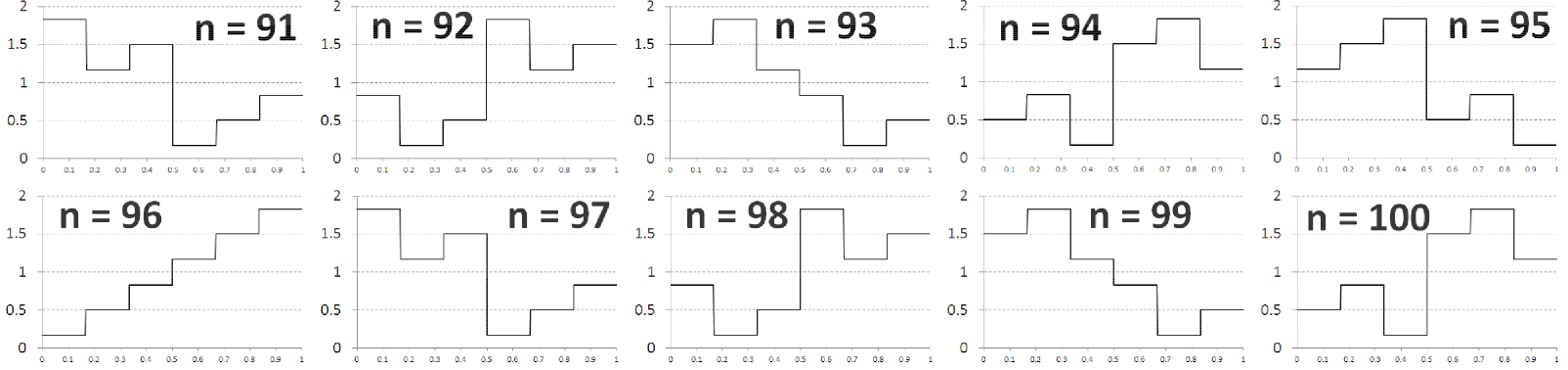}}
\end{center}
\caption
{Results of $((\mathcal{L}_{R^{\tau}})^n f_0)(\cdot )$
for $\tau $ $=(4,5,6,2,3,1)$ and $91 \leq n \leq 100$. }
\label{fig:R_tau3_456231}
\end{figure}
\begin{figure}[htp]
\begin{center}
\rotatebox{-90}{\includegraphics[width=47mm]{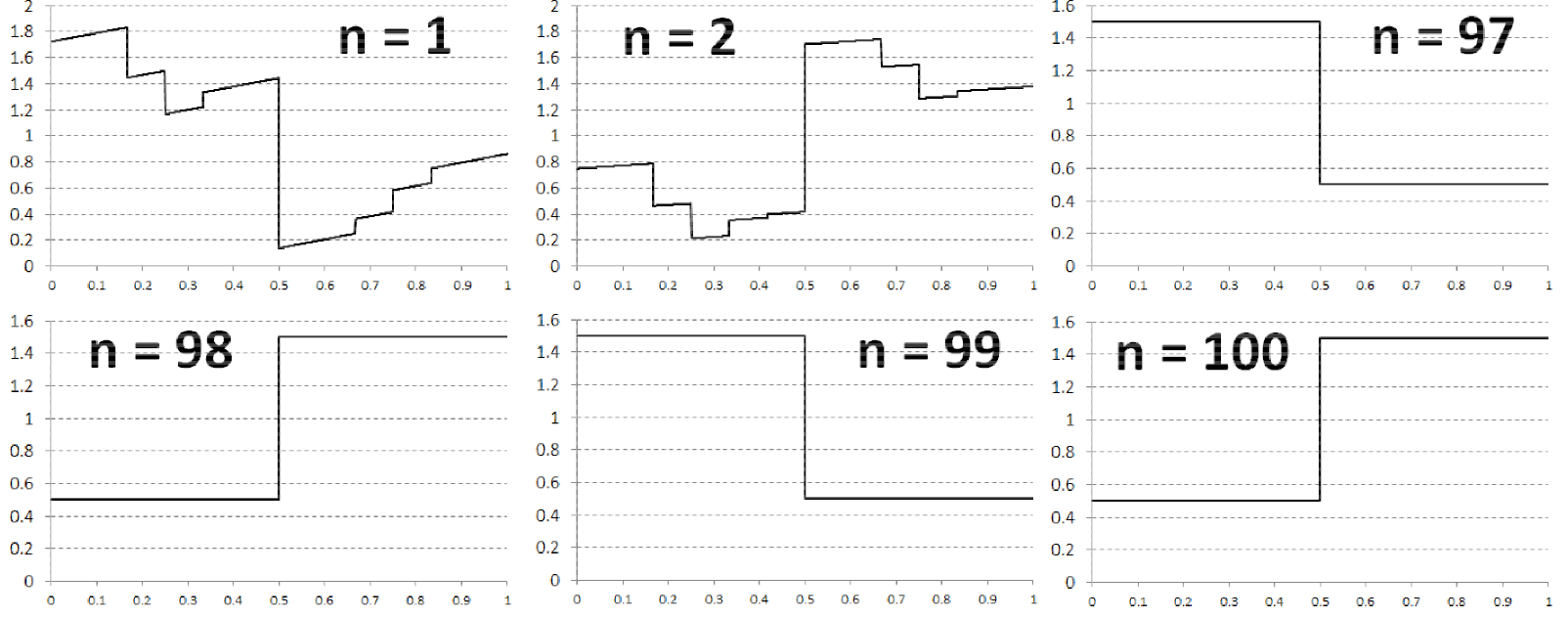}}
\end{center}
\caption
{Results of $((\mathcal{L}_S)^n f_0)(\cdot )$, 
with $T_{y_1}$ $=R^{(3,4,2,1)}$
and $T_{y_2}$ $=R^{(4,5,6,2,3,1)}$, 
for $n=1$, $2$, $97$, $98$, $99$, and $100$. }
\label{fig:R_tau3_3421or456231}
\end{figure}

Here, we consider the skew product $S$ of the transformations
$T_{y_1}=R^{(3,4,2,1)}$ and $T_{y_2}=R^{(4,5,6,2,3,1)}$.
Then, $((\mathcal{L}_{T_{y}})^n f_0)(\cdot )$ $(y=y_1, y_2)$ have
the property of asymptotic periodicity,
and we obtain the following:
\begin{align*}
&Y_1=\{y_1, y_2\}, \ s(T_{y_1})=s(T_{y_2})=1,
\ r(1,T_{y_1})=4, \ r(1,T_{y_2})=6, \\
&g_{1,j_1,T_{y_1}}(x)=4 \times 1_{\left[\frac{j_1-1}{4}, \frac{j_1}{4}\right]}(x),
\ (1\leq j_1 \leq 4),
\ g_{1,j_2,T_{y_2}}(x)=6 \times 1_{\left[\frac{j_2-1}{6}, \frac{j_2}{6}\right]}(x),
\ (1\leq j_2 \leq 6).
\end{align*}
The graphs of $((\mathcal{L}_{T_{y}})^n f_0)(\cdot )$ $(y=y_1, y_2)$, and
$((\mathcal{L}_S)^n f_0)(\cdot )$ are shown in Figs \ref{fig:R_tau3_3421},
\ref{fig:R_tau3_456231}, and \ref{fig:R_tau3_3421or456231}, respectively.
Further, $((\mathcal{L}_S)^n f_0)(\cdot )$ has
the property of asymptotic periodicity, as follows:
\begin{align*}
&\widehat{s}=1, \ \widehat{r}(1)=2,
\ \widehat{g}_{1,j}(x)=2 \times 1_{\left[\frac{j-1}{2}, \frac{j}{2}\right]}(x),
\ (1\leq j \leq 2).
\end{align*}
Note, therefore, that
\begin{align*}
g_{1,T_{y_1}}(x) = g_{1,T_{y_2}}(x)=\widehat{g}_{1}(x)=1_{\left[ 0,1 \right]}(x)
\end{align*}
holds, and $\widehat{r}(1)=2$ is a divisor of $r(1,T_{y_1})=4$ and $r(1,T_{y_2})=6$,
which corresponds to the results given in Theorem \ref{thm:1}.
Because $f_0(x)=2x$ is a monotone increasing function,
we can expect that $0<\lambda_{1,j,T_y}(f_0) < \lambda_{1,j+1,T_y}(f_0)$ holds
for $j =1, \ldots , r(1,T_y)-1$ $(y \in \{y_1, y_2\})$, and that
$0<\widehat{\lambda}_{1,1}(f_0) < \widehat{\lambda}_{1,2}(f_0)$ holds.
Actually, we can confirm these tendencies
in Fig \ref{fig:R_tau3_3421} $[n = r(1,T_{y_1}) \times 25 =100]$,
Fig \ref{fig:R_tau3_456231} $[n = r(1,T_{y_2}) \times 16 = 96]$,
and Fig \ref{fig:R_tau3_3421or456231} $[n = \widehat{r}(1) \times 50 = 100]$.

\subsubsection*{For $T_{y_1}=R^{(2,1,6,4,3,5)}$ and $T_{y_2}=I_d$}

\begin{figure}[htp]
\begin{center}
\rotatebox{-90}{\includegraphics[width=57mm]{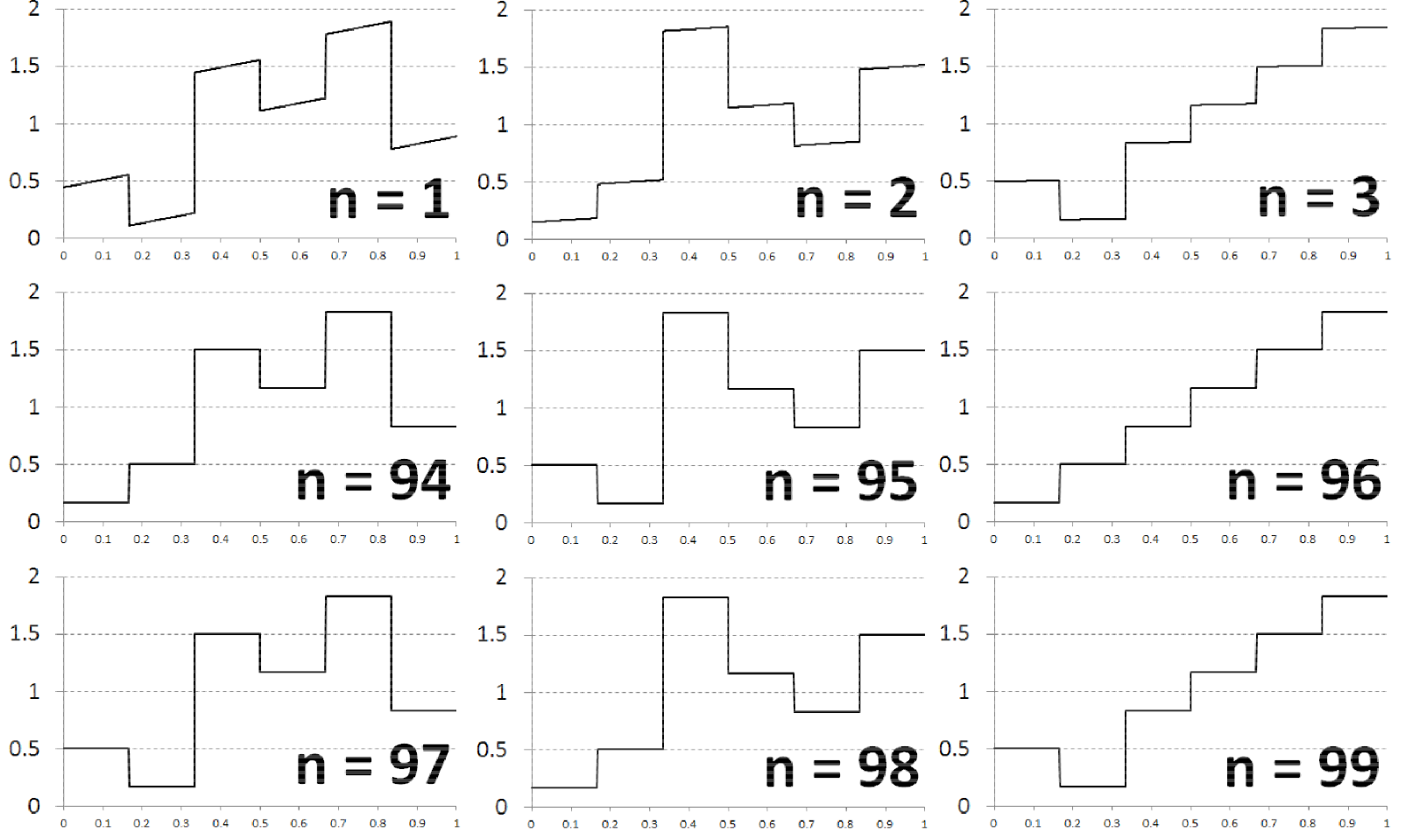}}
\end{center}
\caption
{Results of $((\mathcal{L}_{R^{\tau}})^n f_0)(\cdot )$
for $\tau $= $(2,1,6,4,3,5)$, $n=1, 2, 3, 94, 95, 96, 97, 98$ and $99$. }
\label{fig:LS_R_tau3=(2,1,6,4,3,5)}
\end{figure}
\begin{figure}[htp]
\begin{center}
\rotatebox{-90}{\includegraphics[width=41mm]{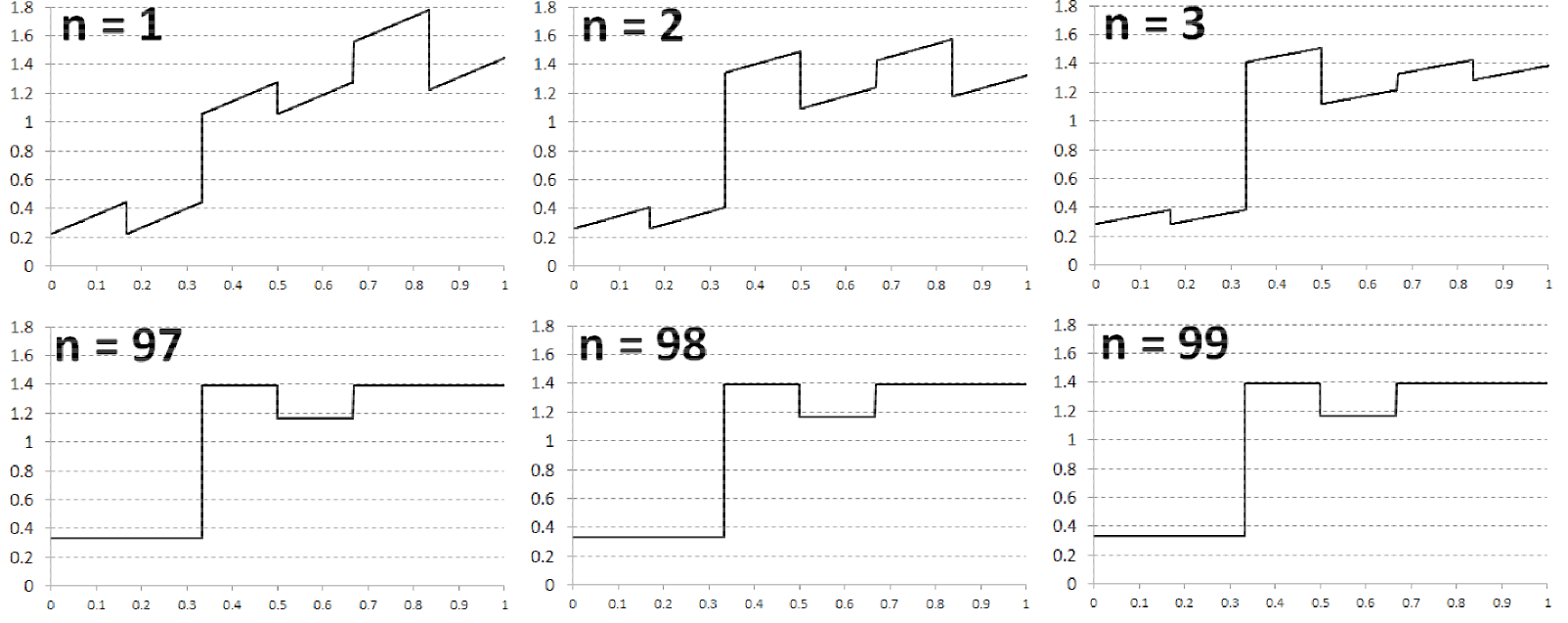}}
\end{center}
\caption
{Results of $((\mathcal{L}_S)^n f_0)(\cdot )$,
with $T_{y_1}$ $=R^{(2,1,6,4,3,5)}$ and $T_{y_2}$ $=I_d$,
for $n=1$, $2$, $3$, $97$, $98$, and $99$. }
\label{fig:LS_R_tau3=(2,1,6,4,3,5)_or_Id}
\end{figure}

Here, we consider the skew product $S$ of the transformations
$T_{y_1}=R^{(2,1,6,4,3,5)}$ and $T_{y_2}=I_d$.
Then, the Perron--Frobenius operator $\mathcal{L}_S f$ is obtained as follows:
\begin{align*}
&(\mathcal{L}_S f)(x)=
\frac{1}{2} (\mathcal{L}_{R^{(2,1,6,4,3,5)}} f)(x) + \frac{1}{2}f(x).
\end{align*}
Because $I_d(x)=x$ $(x\in [0,1])$ is not expanding,
$\mathcal{L}_{I_d}$ does not have the property of asymptotic periodicity,
and $Y_1=\{y_1\}$ holds.
The graphs of $((\mathcal{L}_{R^{(2,1,6,4,3,5)}})^n f_0)(\cdot )$ and
$((\mathcal{L}_S)^n f_0)(\cdot )$ are shown in
Figs \ref{fig:LS_R_tau3=(2,1,6,4,3,5)} and
\ref{fig:LS_R_tau3=(2,1,6,4,3,5)_or_Id}.
Then, we can confirm the following result:
\begin{align*}\allowdisplaybreaks[4]
&s(T_{y_1})=\widehat{s}=3,\ r(1,T_{y_1})=2,\ r(2,T_{y_1})=3,\ r(3,T_{y_1})=1,
\ \widehat{r}(1)=1,\ \widehat{r}(2)=1,\ \widehat{r}(3)=1, \\
&g_{1,1,T_{y_1}}(x)=6 \times 1_{\left[0, \frac{1}{6}\right]}(x),
\ g_{1,2,T_{y_1}}(x)=6 \times 1_{\left[\frac{1}{6}, \frac{2}{6}\right]}(x), \\
&g_{2,1,T_{y_1}}(x)=6 \times 1_{\left[\frac{2}{6}, \frac{3}{6}\right]}(x),
\ g_{2,2,T_{y_1}}(x)=6 \times 1_{\left[\frac{4}{6}, \frac{5}{6}\right]}(x),
\ g_{2,3,T_{y_1}}(x)=6 \times 1_{\left[\frac{5}{6}, 1\right]}(x), \\
&g_{3,1,T_{y_1}}(x)=6 \times 1_{\left[\frac{3}{6}, \frac{4}{6}\right]}(x), \\
&\widehat{g}_{1,1}(x)=3 \times 1_{\left[0, \frac{1}{3}\right]}(x),
\ \widehat{g}_{2,1}(x)=2 \times
1_{\left[\frac{2}{6}, \frac{3}{6}\right]\cup \left[\frac{4}{6}, 1\right]}(x),
\ \widehat{g}_{3,1}(x)=6 \times 1_{\left[\frac{3}{6}, \frac{4}{6}\right]}(x).
\end{align*}
Note, therefore, that
\begin{align*}
& g_{1,T_{y_1}}(x)=\widehat{g}_{1}(x)=3 \times 1_{\left[0, \frac{1}{3}\right]}(x), \\
& g_{2,T_{y_1}}(x)=\widehat{g}_{2}(x)=2 \times
1_{\left[\frac{2}{6}, \frac{3}{6}\right]\cup \left[\frac{4}{6}, 1\right]}(x), \\
& g_{3,T_{y_1}}(x)=\widehat{g}_{3}(x)=6 \times 1_{\left[\frac{3}{6}, \frac{4}{6}\right]}(x)
\end{align*}
hold, and that $\widehat{r}(i)=1$ $(i=1,2,3)$
corresponds to the result given in Theorem \ref{thm:2}.
\end{ex}

\begin{ex}
For a constant $a\in (0, 1/3)$, we define the disjoint subintervals
$I_k$ $(1\leq k \leq 9)$ as
\begin{align*}
&I_k \equiv [c_{k-1}, c_k) \ (1\leq k \leq 8)
\ \mbox{and}\ I_9 \equiv [c_8, c_9], \\
&\quad \mbox{where }c_0=0,\ c_k=\left\{
\begin{array}{ll}
c_{k-1}+a, & (k=1, 5, 9),\\
c_{k-1}+\frac{1-3a}{6}, & (k=2, 3, 4, 6, 7, 8).
\end{array}
\right.
\end{align*}
Note that $X=[0,1]=\bigcup _{k=1}^9 I_k$ holds.
Then, we define the transformation $Q^{(a)}: X \rightarrow X$ as
\begin{align*}
&Q^{(a)} x =\left\{
\begin{array}{ll}
c_1+D_{k-1,k}^{1,4}(x - c_{k-1}), & x \in I_k \ (1\leq k \leq 4),\\
c_1+D_{k-1,k}^{1,8}(x - c_{k-1}), & x \in I_k \ (k=5),\\
c_5+D_{k-1,k}^{5,8}(x - c_{k-1}), & x \in I_k \ (6\leq k \leq 9),
\end{array}
\right. \quad \mbox{where }D_{k_1, k_2}^{k_3, k_4}
= \frac{ c_{k_4} - c_{k_3} }{ c_{k_2} - c_{k_1} }.
\end{align*}
The graph of $\{Q^{(a)} x;\ x\in [0, 1]\}$ for $a=0.05$ is shown in Fig \ref{fig:Graph_Q_a}.
Further, we define the transformation $Q^{(a,b)}$ as
$Q^{(a,b)} x = Q^{(a)}x+b$ $(x \in X)$ for $b \in [-a, a]$.
\begin{figure}[htp]
\begin{center}
\rotatebox{-90}
{\includegraphics[width=48mm]{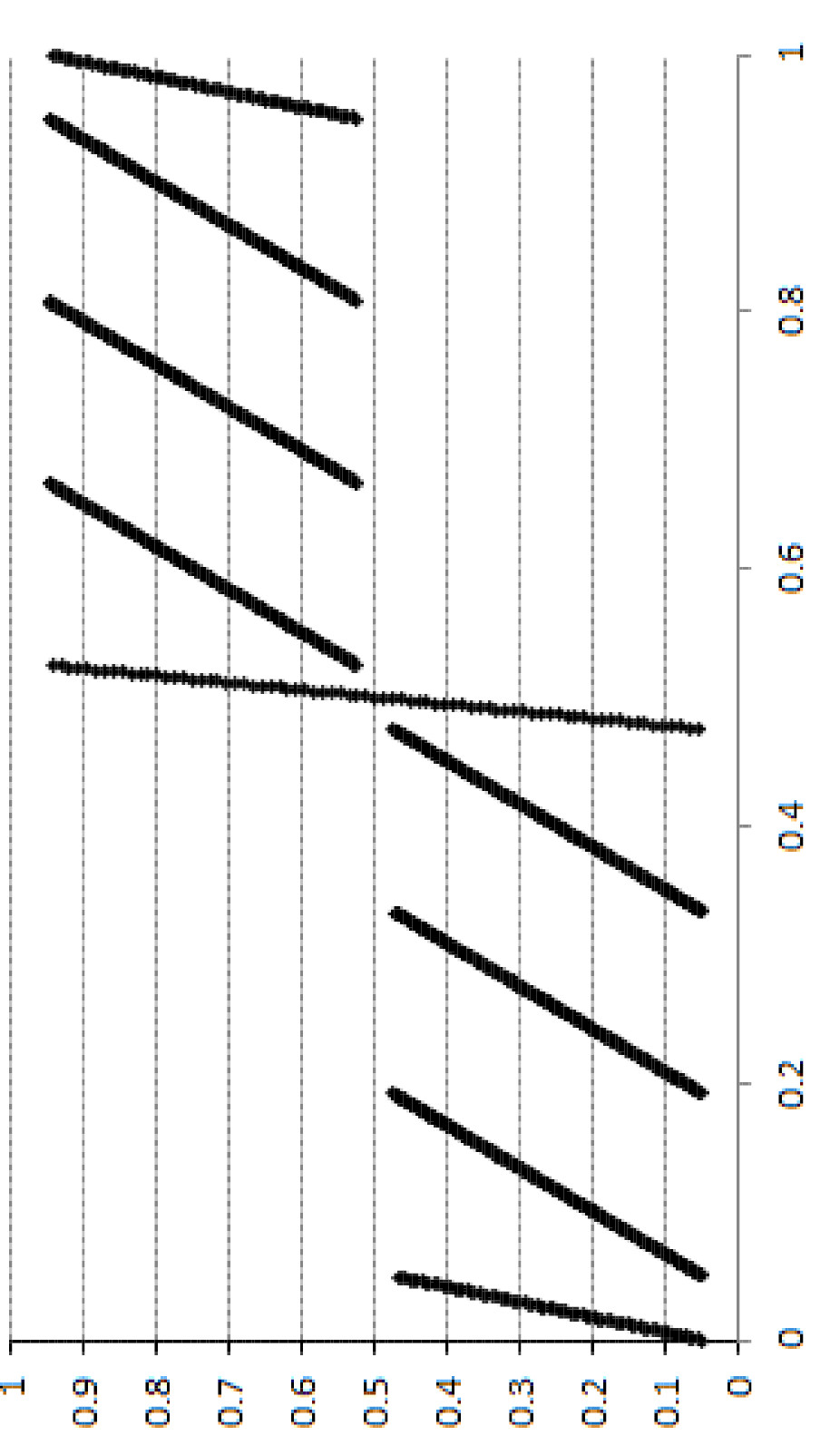}}
\end{center}
\caption
{$\{Q^{(a)}(x);x\in [0, 1]\}$ for $a=0.05$. }
\label{fig:Graph_Q_a}
\end{figure}
Then, $Q^{(a,b)}$ is considered to be the perturbed transformation of $Q^{(a)}$,
and the Perron--Frobenius operator ${\cal L}_{Q^{(a,b)} } f$ is obtained as follows [(i)--(v)]:
\noindent
\begin{align*}\allowdisplaybreaks[4]
\mbox{(i)}&\ ({\cal L}_{Q^{(a,b)}} f)(x) = D_{1,4}^{0,1} f(c_0)\ \mbox{for}\ x = c_1+b; \\
\mbox{(ii)}&\ ({\cal L}_{Q^{(a,b)}} f)(x) = D_{1,8}^{4,5} f(c_4+D_{1,8}^{4,5}(x-b-c_1))
+\sum_{k=1}^4 D_{1,4}^{k-1,k} f\big( c_{k-1}+ D_{1,4}^{k-1,k} (x - b - c_1) \big) \\
&\ \mbox{for}\ x \in (c_1+b, c_4+b]; \\
\mbox{(iii)}&\ ({\cal L}_{Q^{(a,b)}} f)(x) = D_{1,8}^{4,5} f(c_4+D_{1,8}^{4,5}(x-b-c_1))
\ \mbox{for}\ x \in (c_4+b, c_5+b]; \\
\mbox{(iv)}&\ ({\cal L}_{Q^{(a,b)}} f)(x) = D_{1,8}^{4,5} f(c_4+D_{1,8}^{4,5}(x-b-c_1))
+\sum_{k=6}^9 D_{5,8}^{k-1,k} f\big( c_{k-1}+ D_{5,8}^{k-1,k} (x - b - c_5) \big) \\
&\ \mbox{for}\ x \in (c_5+b, c_8+b]; \\
\mbox{(v)}&\ ({\cal L}_{Q^{(a,b)}} f)(x) = 0\ \mbox{for}\ x \in [0, c_1+b)\cup (c_8+b ,1].
\end{align*}

Setting $\displaystyle \widehat{a}\equiv \frac{a(1-3a)}{2(1-2a)}$,
we then obtain the following properties:
\begin{align*}\allowdisplaybreaks[4]
\mbox{(A)}\ \ ( (Q^{(a, b)})^n x_0)_{n=0}^{\infty}
& \subset [0, c_4+\widehat{a}]\ \ \mbox{for}
\ x_0 \in [0, c_4+\widehat{a}] \ \mbox{and}\ b \leq \widehat{a}; \\
\mbox{(B)}\ \ ( (Q^{(a, b)})^n x_0)_{n=0}^{\infty}
& \subset [c_5-\widehat{a}, 1]\ \ \mbox{for}
\ x_0 \in [c_5-\widehat{a}, 1] \ \mbox{and}\ b \geq -\widehat{a}.
\end{align*}
For $a=0.15$, $b_1=a/4$ $(<\widehat{a})$, and $b_2=-3a/4$ $(<-\widehat{a})$,
we consider the transformations $T_{y_1}=Q^{(a,b_1)}$, $T_{y_2}=Q^{(a,b_2)}$,
and the corresponding skew product transformation $S$.
The graphs of $((\mathcal{L}_{Q^{(a,b_k)}})^n f_0)(\cdot )$, $(k=1,2)$,
and $((\mathcal{L}_S)^n f_0)(\cdot )$ are shown in
Figs \ref{fig:Q_ab_(b=0.25a)_(a=0.15)}, \ref{fig:Q_ab_(b=m0.75a)_(a=0.15)},
and \ref{fig:Q_ab_(b=m0.75a_or_0.25a)_(a=0.15)}, respectively.
Then,
\begin{align*}
&Y_1=\{y_1, y_2\}, \ s(T_{y_1})=2,\ s(T_{y_2})=\widehat{s}=1, \\
&r(1,T_{y_1})=r(2,T_{y_1})=1, \ r(1,T_{y_2})=1, \ \widehat{r}(1)=1
\end{align*}
are obtained.
Further, we can expect, and confirm from
Figs \ref{fig:Q_ab_(b=0.25a)_(a=0.15)}--\ref{fig:Q_ab_(b=m0.75a_or_0.25a)_(a=0.15)} $(n=100)$,
that the densities
$g_{i,T_{y_1}}(x)$, $(i=1,2)$,
$g_{1,T_{y_2}}(x)$, and $\widehat{g}_{1}(x)$
satisfy
\begin{align*}
&\{ g_{1,T_{y_1}}>0 \} \approx [0.1873, 0.4626], \quad \{ g_{2,T_{y_1}}>0 \} \approx [0.6126, 0.8873], \\
&\{ g_{1,T_{y_2}}>0 \} \approx [0.0373, 0.3126], \quad \{ \widehat{g}_{1}>0 \} \approx [0.0373, 0.4626] .
\end{align*}
Thus, we have that
\begin{align*}
&\{ g_{1,T_{y_1}}>0 \}
\subsetneqq \{ \widehat{g}_{1}>0 \},
\ \{ g_{2,T_{y_1}}>0 \} \cap \{ \widehat{g}_{1}>0 \} = \emptyset ,
\ \mbox{and}\ \{ g_{1,T_{y_2}}>0 \}
\subsetneqq \{ \widehat{g}_{1}>0 \}
\end{align*}
hold, which correspond to the result given in Theorem \ref{thm:1} (1).
Note also that $Q^{(a,b_2)}$ satisfies property (A),
and $Q^{(a,b_2)}$ does not satisfy property (B).
If $x_0$ $\in$ $[c_5-\widehat{a}, 1]$ satisfies
$(Q^{(a, b_2)})^{n_0^*} x_0\in [0, c_4+\widehat{a}]$ for some $n_0^* \in \mathbb{N}$,
then we have that $( (Q^{(a, b_2)})^n x_0)_{n=n_0^*}^{\infty} \subset [0, c_4+\widehat{a}]$.
Thus, we can expect, and confirm from Fig \ref{fig:Q_ab_(b=m0.75a)_(a=0.15)},
that $\lim_{n\to \infty} (({\cal L}_{Q^{(a,b_2)}})^n f_0)(x)=0$
holds for $x \in [c_5-\widehat{a}, 1]$.

\begin{figure}[htb]
\begin{center}
\rotatebox{-90}
{\includegraphics[width=23.1mm]{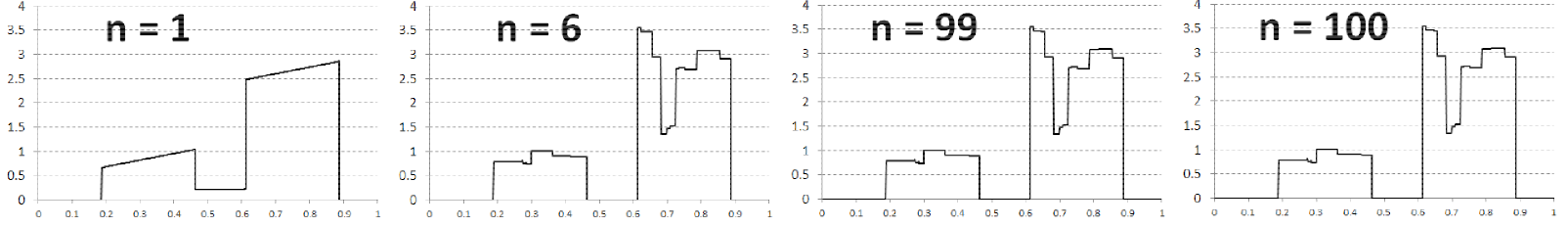}}
\end{center}
\caption
{Results of $((\mathcal{L}_{Q^{(a,b_1)}})^n f_0)(\cdot )$
with $a=0.15$ and $b_1=a/4$\ $(n=1, 6, 99$, and $100)$. }
\label{fig:Q_ab_(b=0.25a)_(a=0.15)}
\end{figure}
\begin{figure}[htp]
\begin{center}
\rotatebox{-90}
{\includegraphics[width=23.1mm]{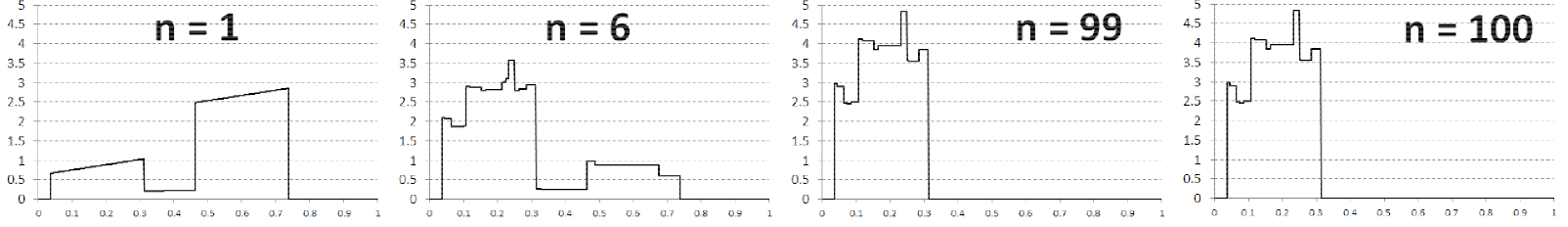}}
\end{center}
\caption
{Results of $((\mathcal{L}_{Q^{(a,b_2)}})^n f_0)(\cdot )$
with $a=0.15$ and $b_2=-3a/4$\ $(n=1, 6, 99$, and $100)$. }
\label{fig:Q_ab_(b=m0.75a)_(a=0.15)}
\end{figure}
\begin{figure}[htp]
\begin{center}
\rotatebox{-90}
{\includegraphics[width=23.1mm]{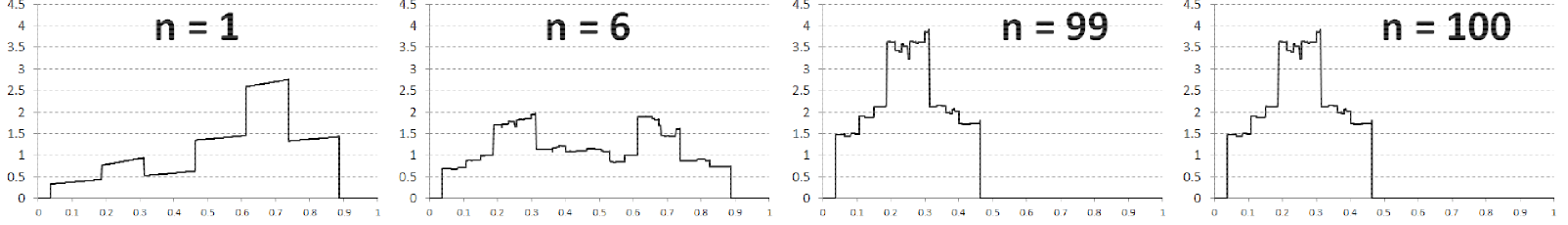}}
\end{center}
\caption
{Results of $((\mathcal{L}_S)^n f_0)(\cdot )$
for $T_{y_k}=Q^{(a,b_k)}$ $(k=1,2)$
with $a=0.15$, $b_1$ $=a/4$, and $b_2$ $=-3a/4$
$(n=1, 6, 99$, and $100)$. }
\label{fig:Q_ab_(b=m0.75a_or_0.25a)_(a=0.15)}
\end{figure}
\end{ex}

\section{Conclusions}

In this paper, we studied the effects of randomization on the asymptotic periodicity. We showed that the supports of the ergodic probability densities for random iterations include at least one support of the ergodic probability density for almost all $m$-nonsingular transformations. This implies that the number of ergodic components of random iterations is not greater than the number of ergodic components of each of the $m$-nonsingular transformations. We also discussed the period of the limiting densities of random iterations. Our results suggest that even a small noise could change the ergodic properties of the system.

\section*{Acknowledgments}
The authors thank the reviewer for
various comments and constructive suggestions to improve the quality of the paper.

\begin{flushleft}
\mbox{  }\\
\hspace{61mm} Hiroshi Ishitani\\
\hspace{61mm} Department of Mathematics\\
\hspace{61mm} Mie University\\
\hspace{61mm} Tsu, Mie 514-8507\\
\hspace{61mm} Japan\\
\hspace{61mm} e-mail: ishitani@edu.mie-u.ac.jp\\
\mbox{  }\\
\hspace{61mm} Kensuke Ishitani\\
\hspace{61mm} Department of Mathematics and Information Sciences\\
\hspace{61mm} Tokyo Metropolitan University\\
\hspace{61mm} Hachioji, Tokyo 192-0397\\
\hspace{61mm} Japan\\
\hspace{61mm} e-mail: k-ishitani@tmu.ac.jp
\end{flushleft}

\end{document}